\documentclass[11pt]{article}
\usepackage{amsmath,amsthm,amssymb}
\usepackage[mathscr]{eucal}
\usepackage[utf8]{inputenc}
\usepackage{graphicx}
\usepackage{bm}
\usepackage{subcaption}
\usepackage{color}

\usepackage[english]{babel}
\usepackage{amsmath,amstext,amssymb}
\usepackage{array}

\title{A general approach to optimal imperfect maintenance
activities of a repairable equipment with imperfect
maintenance and multiple failure modes}
\author{Rubén Mullor\\Dpto. Matemáticas, Universidad de Alicante\\ \texttt{ruben.mullor@ua.es}
\and Julio Mulero \\Dpto. Matemáticas, Universidad de Alicante\\ \texttt{julio.mulero@ua.es}
\and Mario Trottini \\Dpto. Matemáticas, Universidad de Alicante\\ 
\texttt{mario.trottini@ua.es}
}
\date{}

\numberwithin{equation}{section} \numberwithin{theorem}{section}

\usepackage{eurosym} 

\begin{document}

\maketitle
\begin{abstract} 
In this paper we describe  a general approach to optimal imperfect maintenance activities of a repairable equipment with independent components. Most of the existing works on optimal imperfect maintenance activities of a repairable equipment with independent components. In addition, it is assumed  that all the components of the equipment share the same model and the same maintenance intervals and that effectiveness of maintenance is known. In this paper we take a different approach.  In order to formalize the uncertainty on the occurrence of failures and on the effect of maintenance activities we consider, for each component,  a class of candidate models  obtained combining models for failure rate with models for imperfect maintenance  and let the data select the best model (that might be different for the different components).  All the parameters are assumed to be unknown and are jointly estimated via maximum likelihood. Model selection is performed, separately for each component, using  standard selection criteria that take into account the problem of over-parametrization.  The selected models are used  to derive the cost per unit time and the average reliability of the equipment, the objective functions of a Multi-Objective Optimization Problem with maintenance intervals of each single component as decision variables. The proposed procedure is illustrated using a real data example.

\vspace*{0.5cm}
\noindent \textbf{Keywords} Reliability, Imperfect Maintenance,  Multiple failure modes, Ageing Process, Maximum Likelihood, Multiobjective Optimization.
\end{abstract}
%

\section{Introduction}
\label{sec1}

During the last years, the need of optimizing preventive maintenance  activities, performed in different repairable equipments of industrial plants, has experimented a great increase as a consequence of the interest of improving plants economy and safety (see, for example, Qiu et al., 2017a, or Zheng
et al., 2016). Cost  and reliability of a reparable equipment, in fact, depend on the random occurrences of failures which, in turn, depend on  the age  of the equipment  which  is a function of its chronological time  and  the preventive maintenance activities performed to slow down its deterioration.  For an equipment with $q$ independent components, if preventive maintenance activities for the $k$-th component are scheduled  every $M_{k}$ units of time, (i.e. if maintenance intervals are constant), the goal is to find the optimal value of $(M_{1},\ldots,M_{q})$ that provides a suitable trade-off between the two conflicting objectives of minimize cost and maximize reliability (see Wang and Pham,
2006 for a complete revision on maintenance methods).

A suitable solution of the problem requieres a separate modeling, for each component of the equipment, of (i)  the failure rate; and (ii)  the effect of maintenance activities on the age of the component. 

In an optimal maintenance policy setting,  the \textit{Weibull model} and \textit{imperfect maintenance} (IM) models have played a key role in modeling the failure rate and the effect of maintenance activities on the age of a component. In fact, the most popular failure rate model is the {\em Weibull} one which includes as special case the exponential and linear failure rate model (see Tan
and Kramer, 1995). 

IM models assumes that maintenance activities reduce the equipment age by some degree depending on their effectiveness. The most important IM models are the \textit{Proportional Age Setback} (PAS), proposed by Martorell et al. (1999), and \textit{Proportional Age Reduction} (PAR), proposed by Malik (1979).  While the PAS model considers that each maintenance activity reduces the total equipment age by a factor $\varepsilon$, the PAR model assumes that only the age gained between two consecutive maintenance activities is reduced by a factor $\varepsilon$ (see also Kijima and Sumita, 1986; Kijima,
1989; or Tanwar et al., 2014, for a more recent revision). Both models   include as special, and extreme cases, for a suitable choice of the parameter $\varepsilon$,  the less realistic  \textit{Good as New} (GAN) model in  which the the age of the equipment after each maintenance is restored to the initial time ($\varepsilon=1$); and the \textit{Bad as Old} (denoted as BAO) model, which assumes that the maintenance action has no effect on the age of the equipment ($\varepsilon=0$). 

Most of the existing literature on optimal maintenance policy considers a single model for equipment behavior obtained by combining a model for the random occurrences of failures and a model for the effect of maintenance activities on the age of the component. In addition it is assumed  that all the components of the equipment share the same model and the same maintenance interval (i.e., $M_{1}=M_{2}\ldots=M_{q}$ in the optimization problem). A common choice for the maintenance model is  BAO or GAN (as in Neves
et al., 2011, Taghipour and Banjevic, 2013 or Qiu et al., 2017b). IM models,  although more realistic, are less common and, when used, the parameter that represents the  degree  by which maintenance reduces the age of the equipment is assumed to be known (as in Tsai et al., 2011).

Effectiveness of maintenance, and failure rate,  however are usually unknown and might be different for different components of a repairable equipment. This suggests that: (a) considering a single model for equipment behaviour introduces unnecessary  assumptions on the random occurrence of failures and on the effectiveness of maintenance, that could be avoided by considering a richer class of models; (b) assuming that the effectiveness parameter in IM models is known underestimates the uncertainty of the problem.
For a given model for equipment behaviour, the effectiveness parameter should be considered unknown and estimated simultaneously with the parameter of the failure rate model; (c)  different components should be allowed to follow different models and have different maintenance intervals.

In order to address the limitations of current proposals (see, for instance, Taghipour and Banjevic, 2013), in this paper we propose a general approach to the problem that can be described in terms of the following four step procedure:
\vspace{0.2cm}

\noindent   {\bf STEP 1}  For each equipment component we consider a class of four candidate models {\em PAS-linear}; {\em PAR-linear}, {\em PAS-Weibull}, {\em PAR-Weibull} obtained by combining two possible  models for failure rate, {\em linear} and {\em Weibull}; with two possible models for maintenance, PAR and PAS.   Both the parameters of the failure rate models and  the effectiveness parameter $\varepsilon$ of the IM models are assumed to be unknown and need to be estimated from the data (see
Wang and Yang, 2012, and S\'anchez et al., 2006). Note that despite the fact that the Weibull failure rate model includes, as special case,  the simpler linear failure rate model, when the model's parameter are unknown the linear model should be considered as a separate option. If the ``true'' failure rate is linear, a Weibull model would be over-parametrized, and the estimated linear model would be preferable, for inferential purposes, to the estimated Weibull model (the two models are asimptotically equivalent but might provide very different results if sample size is small). 
\vspace{0.2cm}

\noindent   {\bf STEP 2}   For each model defined in Step 1,  the model's parameters (which include the parameters of the failure rate models and   the effectiveness parameter $\varepsilon$ of the IM models) are estimated simultaneously by  maximum likelihood using a direct search algorithm based on the \textit{Nelder-Mead Simplex} (NMS) method (see, for instance, Nelder and Mead, 1965; Lagarias et al.,
1998; S\'anchez et al., 2006).

\vspace{0.2cm}

\noindent   {\bf STEP 3}   For each component of the equipment, we select the best model in the class of models defined in Step 1. As observed before, most of the existing literature on  optimal maintenance policy considers a single model  for equipment behaviour. As a result the problem of model selection is usually not addressed. 
An important exception is Tan and Kramer (1995) which consider a class of  BAO-GAN models and select the best model using Leave-One-Out cross validation. Here we 
implement Leave-One-Out cross validation for the richer class of models defined in Step 1, and introduce two additional selection criteria de \textit{Akaike Information Criterion} and  the  \textit{Bayesian Information Criterion}  (Hastie et al., 2008) to address the issue of over-parametrization.
\vspace{0.2cm}

\noindent   {\bf STEP 4}   We use the selected models in Step 3  as  input of a  \textit{Multi-Objective Optimization Problem} where average reliability and cost per unit of time are the objective functions, and maintenance intervals of each single component the  decision variables (see S\'anchez et al., 2007; Feng et al., 2017; and Tsai
et al., 2011).  The optimization problem is solved applying an algorithm for nonlinear constrained optimization, based on Sequential Quadratic Programming (SQP) method (see, for example, Alrabghi and Tiwari, 2015, and Xia
et al., 2012).
\vspace{0.2cm}

The rest of the paper is organized as follows.  In Section \ref{sec2} we introduce the  {\em PAS-linear}; {\em PAR-linear}, {\em PAS-Weibull}, {\em PAR-Weibull} models and derive the corresponding expression of hazard and cumulative hazard functions.  In Section \ref{sec3} we  obtain the likelihood function four the four models and discuss maximum likelihood estimation. The problem of over-parametrization and model selection is also addressed in Section \ref{sec3}. 
In Section \ref{sec4} we obtain the expression of the average reliability and   cost per time units for each of the four model and 
describe the optimization algorithm. 

An illustration of the proposed approach is presented in Section \ref{sec5}
using a data set  containing failure time and maintenance activities of eight motor-operated safety valves of a Nuclear Power Plant (NPP) over a period of ten years. Finally,
 Section \ref{sec6} summarizes the main results.

\section{Reliability models under imperfect maintenance}
\label{sec2}

In our modelling of equipment behaviour we consider two possible  models for failure rate, {\em linear} and {\em Weibull}, and  two possible models for maintenance: the \textit{Proportional Age Setback} (PAS) and the \textit{Proportional Age Reduction} (PAR) models. As observed in the introduction, both PAR and PAS models assume that each maintenance activity reduces the age of the equipment by some degree depending  on its effectiveness and provide a much more realistic alternative to the 
extreme ``Bad as Old'' (BAO) and ``Good as New'' (GAN) maintenance models to explain the effect of maintenance activities on the state of equipments. This produces a class of four candidate models: {\em PAS-linear}; {\em PAR-linear}, {\em PAS-Weibull}, {\em PAR-Weibull}.
For each of the four models it is assumed that working conditions are normal and initial hazard rate is null. These are reasonable assumptions for the application that we discuss in section \ref{sec5} and simplify  the expressions involved in parameters estimation. A more detailed description of the four models follows.

\subsection{Proportional Age Setback model}\label{sec2:subsec1}

In the PAS approach, each maintenance activity is assumed to shift the origin of time from which the age of the equipments are evaluated. The PAS model considers that maintenance activities reduce proportionally to a factor $\varepsilon$, the age the equipment had immediately before it enters in maintenance. Following the model characterization, the age of the equipment in instant $t$ of period $m$, $w_{m-1}^+$, after the $(m-1)$-maintenance activity is given by
\begin{equation}\label{eq1}
\omega_{m-1}^+=t-\sum_{k=0}^{m-2}(1-\varepsilon)^{k}\varepsilon\tau_{m-k-1}, 
\end{equation}
where $\tau_{j}$ is the time in which the equipment undertakes the $j$-maintenance activity.

Using \eqref{eq1}, it is possible to obtain an age-dependent reliability model in which the induced or conditional failure rate, in period $m$ after the $(m-1)$-maintenance, is given by
\begin{equation}\label{eq2}
h_m(\omega)=h\left(t-\sum_{k=0}^{m-2}(1-\varepsilon)^{k}\varepsilon\tau_{m-k-1}\right). 
\end{equation}

\begin{description}
\item[PAS-linear model]
\end{description}
Considering the age of the equipment after the $(m-1)$-maintenance given by \eqref{eq1}, and adopting a linear model for the failure rate, 
\begin{equation}\label{eq3}
h(\omega)=\alpha \omega(t,\varepsilon),
\end{equation}
the expression for the induced failure rate after the $(m-1)$-maintenance becomes 
\begin{equation*}\label{eq5}
h_m(\omega)=\alpha\left(t-\sum_{k=0}^{m-2}(1-\varepsilon)^{k}\varepsilon\tau_{m-k-1}\right), 
\end{equation*}
where $\alpha$ is the linear aging rate.

The cumulated hazard function $H_m(\omega)$ in the period $m$, after the $(m-1)$-maintenance activity, can be obtained by integration of the induced hazard function obtaining
\begin{equation*}\label{eq6}
H_m(\omega)=\frac{\alpha}{2}\left(t-\sum_{k=0}^{m-2}(1-\varepsilon)^{k}\varepsilon\tau_{m-k-1}\right)^2 . 
\end{equation*}

\begin{description}
\item[PAS-Weibull model]
\end{description}

If, instead of the linear failure rate, a Weibull failure rate with shape parameter $\beta$ and scale parameter $\eta$, is considered, 
\begin{equation}\label{eq7}
h(\omega)=\frac{\beta}{\eta^\beta}(\omega(t,\varepsilon))^{\beta-1}, 
\end{equation}
then, replacing the expression corresponding to $\omega(t,\varepsilon)$ given by \eqref{eq1} into \eqref{eq7}, the expressions for the induced failure rate, and the cumulated hazard function in period $m$ after the $(m-1)$-maintenance  become, respectively, 
\begin{eqnarray*}
h_m(\omega)&=&\frac{\beta}{\eta^\beta}\left(t-\sum_{k=0}^{m-2}(1-\varepsilon)^{k}\varepsilon\tau_{m-k-1}\right)^{\beta-1},  \\
H_m(\omega)&=&\frac{1}{\eta^\beta}\left(t-\sum_{k=0}^{m-2}(1-\varepsilon)^{k}\varepsilon\tau_{m-k-1}\right)^\beta. \\ 
\nonumber
\end{eqnarray*}

\subsection{Proportional Age Reduction model}\label{sec2:subsec2}

In the PAR approach, each maintenance activity is assumed to reduce, proportionally to its effectiveness, the age gained from the previous maintenance. Thus, while the PAS model considers that each maintenance activity reduces the total equipment age, the PAR model assumes that maintenance reduces the age gained between two consecutive maintenance activities by a factor $\varepsilon$.
 
According with the above conditions, the age of the equipment in instant $t$ of period $m$, after the $(m-1)$-maintenance activity using the PAR model is given by
\begin{equation}\label{eq10}
\omega_{m-1}^+=t-\varepsilon\tau_{m-1}. 
\end{equation}

Using a similar procedure as the one described for the PAS model, but adopting \eqref{eq10} instead of \eqref{eq1}, it is possible to derive the expression for the failure rate and the cumulated hazard function of imperfect maintenance in instant $t$, under the PAR approach assuming either the linear or Weibull failure rate model, this produces the PAR-linear and PAR-Weibull models:
\begin{description}
\item[PAR-linear model]
\end{description}
\[
h_m(t)=\alpha(t-\varepsilon\tau_{m-1}); \,\,\,\,\,\,\,\,\,\,\,\, H_m(t)=\frac{\alpha}{2}(t-\varepsilon\tau_{m-1})^2.  
\]
\begin{description}
\item[PAR-Weibull model ]
\end{description}
\[
h_m(t)=\frac{\beta}{\eta^\beta}(t-\varepsilon\tau_{m-1})^{\beta-1}; \,\,\,\,\,\,\,\,\,\,\,\, H_m(t)=\frac{1}{\eta^\beta}(t-\varepsilon\tau_{m-1})^\beta.
\]

\section{Maximum Likelihood Estimation}\label{sec3}

For each of the  four  models for equipment behaviour  discussed in the previous section, estimation of model's parameters is obtained using the Maximum Likelihood procedure (MLE). This  requires an explicit formula for the likelihhood of each model and an optimization algorithm to maximize the likelihood function.  

The likelihood function ($L$) is the product of probabilities of the observed data as a function of the model parameters. For a reliability model, as described in Tan and Kramer (1995), it can be written as
\begin{equation}\label{eq16}
L(\xi|\text{observed data})=\prod_{\textup{failures}}h(t)\prod_{\textup{maint.}}R(t),
\end{equation}
where $\xi$ denotes de vector of model's parameters (i.e. $\xi=(\alpha,\varepsilon)$  for  linear failure rate and  $\xi=(\beta,\eta,\varepsilon)$ for  Weibull  failure rate)
and $R(t)$ is the reliability function, $R(t)=e^{-H(t)}$. 

In particular, for $P$ equipments with a single component under imperfect preventive maintenance,  (\ref{eq16}) becomes
\begin{equation*}\label{eq17}
L(\xi)=\prod_{p=1}^P \left[\prod_{m=1}^{M_p+1}\left\{\prod_{j=1}^{r_{p,m}}h_{p,m}(t_{p,m,j})\exp\left(-\sum_{m=1}^{M_p}H_{p,m}(\tau_{p,m})-H_{M_p+1}(\tau^*_p)\right)\right\}\right],
\end{equation*}
where for each component $p$: $r_{p,m}$  denotes the number of failures of the component during the $m$-maintenance which occur at times $t_{p,m_1},t_{p,m_2},\dots$; $\tau_{p,m}$ is the  chronological time for the $m$-maintenance;   
$M_p$ denotes  the number of preventive maintenance activities performed during the observation period $\tau_p^*$;  $h_{p,m}(t)$ and $H_{p,m}(\tau)$ are the induced hazard function and the cumulated hazard function in period $m$, respectively;  and $H_{M_p+1}(\tau_p^*)$ is the cumulated hazard function in censoring time $\tau_p^*$.

As usual in this context, and for computational purpose, it is preferable to maximize the logarithm of the likelihood function given by
{\small
\begin{equation}\label{eq18}
\log L(\xi)=\sum_{p=1}^P \left[\sum_{m=1}^{M_p+1}\left\{\sum_{j=1}^{r_{p,m}}\log h_{p,m}(t_{p,m,j})-\sum_{m=1}^{M_p}H_{p,m}(\tau_{p,m})-H_{M_p+1}(\tau^*_p)\right\}\right].
\end{equation}}

Replacing $h_{p,m}(t)$ and $H_{p,m}(t)$  in  \eqref{eq18} with the corresponding expressions obtained in Section \ref{sec2}   
produces the following log likelihood functions for the PAS-linear, the PAR-linear, the PAS-Weibull and the PAR-Weibull models:

\medskip
\noindent\textbf{PAS-linear model:}
{\small
\begin{eqnarray*}
\log L(\alpha,\varepsilon)&=&\sum_{p=1}^P \left\{\sum_{m=1}^{M_p+1}\left[\sum_{j=1}^{r_{p,m}}\log \left[\alpha\left(t_{p,m,j}-\sum_{k=0}^{m-2}(1-\varepsilon)^k\varepsilon\tau_{p,m-k-1}\right)\right]\right.\right.\\
&-&\frac{\alpha}{2}\left[\sum_{m=1}^{M_p}\left(\tau_{p,m}-\sum_{k=0}^{m-2}(1-\varepsilon)^k\varepsilon\tau_{p,m-k-1}\right)^2\right.\\
&+&\left.\left.\left.\left(\tau_{p}^*-\sum_{k=0}^{M_p-2}(1-\varepsilon)^k\varepsilon\tau_{p,m-k-1}\right)^2\right]\right]\right\}.
\end{eqnarray*}}

\medskip
\noindent\textbf{PAR-linear  model:}
{\small
\begin{eqnarray*}
 \log L(\alpha,\varepsilon)&=&\sum_{p=1}^P \left\{\sum_{m=1}^{M_p+1}\left[\sum_{j=1}^{r_{p,m}}\log [\alpha(t_{p,m,j}-\varepsilon\tau_{p,m-1})]\right.\right.\\
&-&\left.\left.\frac{\alpha}{2}\left[\sum_{m=1}^{M_p}(\tau_{p,m}-\varepsilon\tau_{p,m-1})^2+(\tau_{p}^*-\varepsilon\tau_{p,m-1})^2\right]\right]\right\}.
\end{eqnarray*}}

\medskip
\noindent\textbf{PAS-Weibull  model:}
\hspace{-0.5cm}{\footnotesize
\begin{eqnarray*}
\log L(\beta,\eta,\varepsilon)&=&\sum_{p=1}^P \left\{\sum_{m=1}^{M_p+1}\left[\sum_{j=1}^{r_{p,m}}\log \left[\frac{\beta}{\eta^\beta}\left(t_{p,m,j}-\sum_{k=0}^{m-2}(1-\varepsilon)^k\varepsilon\tau_{p,m-k-1}\right)^{\beta-1}\right]\right.\right.\\
&-&\frac{1}{\eta^\beta}\left[\sum_{m=1}^{M_p}\left(\tau_{p,m}-\sum_{k=0}^{m-2}(1-\varepsilon)^k\varepsilon\tau_{p,m-k-1}\right)^\beta\right.\\
 &+&\left.\left.\left.\left(\tau_{p}^*-\sum_{k=0}^{M_p-2}(1-\varepsilon)^k\varepsilon\tau_{p,m-k-1}\right)^\beta\right]\right]\right\}.
\end{eqnarray*}}

\medskip
\noindent\textbf{PAR-Weibull  model:}
\hspace{-0.2cm}{\small
\begin{eqnarray*}
\log L(\beta,\eta,\varepsilon)&=&\sum_{p=1}^P \left\{\sum_{m=1}^{M_p+1}\left[\sum_{j=1}^{r_{p,m}}\log \left[\frac{\beta}{\eta^\beta}(t_{p,m,j}-\varepsilon\tau_{p,m-1})^{\beta-1}\right]\right.\right.\\
&-&\left.\left.\frac{1}{\eta^\beta}\left[\sum_{m=1}^{M_p}(\tau_{p,m}-\varepsilon\tau_{p,m-1})^\beta+(\tau_{p}^*-\varepsilon\tau_{p,m-1})^\beta\right]\right]\right\}.
\end{eqnarray*}}

MLE estimations of  the vector of parameters for each of the four model are  obtained by maximizing  the corresponding log likelihood function (as a function of the vector of parameters). To perform the maximization, in this paper we use the Nelder Mead Simplex method (for further details, see
Nelder and Mead, 1965 and Lagarias et al., 1998).

\subsection{Model selection}\label{sec3:subsec2}
Each of the four reliability models   discussed in the previous sections provides a description of the equipment behaviour that in turns can be used to optimize the maintenace intervals. Different models, however, provide different descriptions and thus some model selection procedure is needed before 
optimization of maintanace interval is performed. The criterion of comparing the maximized likelihood, sometimes used in literature, works fine if the competing models have the same number of parameters, but might lead to overparametrization in other cases. It is well known that more complex models tend to have higher likelihood. In this paper, we consider three standard selection criteria that address the issue of overparametrization: {\em Leave One Out Cross Validation}, 
the {\em Akaike information criterion}    and the {\em Bayesian information criterion}.

Leave-one-out Cross-Validation (LCV) is a general technique used to measure how well a certain model generalizes to new data. Ideally, if we had  a large data set, we  could partition the data  in two parts. Use the first part as learning sample  to estimate the different models, and the second part as validation data set to estimate the models fit to new data and select as best model the one with best fit in the validation data set.   However since data set for estimation of reliability models are often  small, this is usually not possible for the applications that we are interested here. To overcome this problem  in  LCV  each data point is considered in turn as validation sample and the rest of the data as learning sample.
 In this paper we use the implementation of LCV proposed by Tan and Kramer (1995).

A drawback of cross validation is that it is computationally intensive. Two simpler alternatives are the  Akaike information criterion (AIC)   and the  Bayesian information criterion (BIC) introduced respectively by Akaike
(1974) and Schwarz (1978). Both criteria provide a tradeoff
 between model complexity and how well the model fits the data. Both criteria can be expressed as the sum of a godness-of-fit term (which is expressed as minus twice the negative of the maximized loglikelihood) and a penalty term which is an increasing function of the number of parameters in the model. 
AIC can be seen as an asymptotic aproximation of the LCV (much simpler to compute) while BIC has a Bayesian interpretation in terms of posterior probabilities of the model (for a detailed description of the two criteria, their properties and their relationship with LCV see, for example, chapter 7 of Hastie et al.,
2008). 

\section{Optimization process}\label{sec4}

As a first step in the optimization process, for each of the four models described in section 2 (PAS-linear, PAS-Weibull, PAR-linear, PAR-Weibull), we obtain the expression of the average hazard,  average reliability and cost in  terms of the maintenance interval
of an equipment with a single component. The maintenance interval, $M$,  is assumed to be constant. 
The results presented easily generalize to the case of equipments  with $q\geq 2$ independent components. In this case, in fact,  the average reliability  of the equipment  is just the  product  of the average reliability  of the single components and the hourly cost the sum of the hourly cost of the single components.  The methodology to construct the average hazard and reliability functions is slightly different depending on the maintenance model used, so we distinguish the two cases.

\subsection{Average hazard and reliability functions for PAS models} \label{AVhPAS}

 As shown in  Martorell et al. (1999) for the PAS model   after few maintenances  the  age of the equipment (as a function of the chronological time) exhibits a stationary behavior. In particular, 
 assuming a constant maintenance interval $M$, the time in which the equipment undertakes the $m$-maintenance activity is $mM$ and the age of the equipment in instant $t$ of period $m$, $w_{m-1}^+$, after the $(m-1)$-maintenance in equation (\ref{eq1}) becomes: 
\begin{equation}\label{eq1BIS}
\omega_{m-1}^+=t-\sum_{k=0}^{m-2}(1-\varepsilon)^{k}\varepsilon  (m-k-1)M, 
\end{equation}
and, for $m$ large enough, (\ref{eq1BIS}) can be written as: 
\begin{equation}\label{eqt26}
\omega_{m-1}^+=t-mM+\frac{M}{\varepsilon}.
\end{equation}
The stationary behaviour of the age of the equipment in the PAS model is illustrated in  Figure 1-(a). 
As a consequence of this stationarity property, average hazard and reliability for the PAS model can be computed by integrating the hazard and reliability function over an arbitrary maintenance interval. 
Substituting age by its expression given by \eqref{eqt26}  in formulas \eqref{eq3} and  \eqref{eq7}, the hazard function for the the PAS-linear and PAS-Weibull models in period $m$ becomes
\begin{equation}\label{eqt29}
h_{m}(t)=\left\{\begin{array}{ll}
\alpha\left(t-mM+\frac{M}{\varepsilon}\right), &\text{for PAS-linear},\\
\frac{\beta}{\eta^\beta}\left(t-mM+\frac{M}{\varepsilon}\right)^{\beta-1}, &\text{for PAS-Weibull}.\\
\end{array}\right.
\end{equation}
The  corresponding averaged hazard functions, $h^*$,  are obtained  integrating  \eqref{eqt29} over an arbitrary maintenance interval and dividing by $M$, the length of the maintenance interval,
\begin{equation*}
h^*(M)=
\frac{1}{M}\int_{(m-1)M}^{mM}h_{m}(t)\ dt=\left\{\begin{array}{ll}
\frac{M\alpha(2-\varepsilon)}{2\varepsilon}, &\text{ for PAS-linear},\\
\frac{M^{\beta-1}}{(\varepsilon\eta)^\beta}\left(1-(1-\varepsilon)^\beta\right), &\text{ for PAS-Weibull}.\\
\end{array}\right.\\
\end{equation*}
Similarly, given the general expression of the reliability function for the linear and Weibull failure models
\begin{equation}\label{eqt39}
R(\omega)=\left\{\begin{array}{ll}
\exp\left(-\frac{\alpha}{2}\omega^2\right), &\text{ for a linear mdel},\\
\exp\left(-\left(\frac{\omega}{\eta}\right)^\beta\right), &\text{ for a Weibull model},
\end{array}\right.
\end{equation}
substituting in (\ref{eqt39})  age by its expression given by \eqref{eqt26}, the  reliability function, for the PAS-linear and PAS-Weibull models, in period $m$ becomes
\begin{equation}\label{eqt41}
R_{m}(t)=\left\{\begin{array}{ll}
\exp\left(-\frac{\alpha}{2}\left(t-mM+\frac{M}{\varepsilon}\right)^2\right), &\text{for PAS-linear}.\\
\exp\left(-\left(\frac{t-mM+\frac{M}{\varepsilon}}{\eta}\right)^\beta\right), &\text{for PAS-Weibull},
\end{array}\right.
\end{equation}
and, as a consequence of stationarity in the PAS model,  the averaged reliability functions, $R^*$, for the two  models is obtained by integrating \eqref{eqt41} in the maintenance interval $((m-1)M,mM)$.

Note that  since the expressions in \eqref{eqt41} are not integrable functions we can use their power series representation to perform the integration. The corresponding approximation of the average reliability $R^*(M)$ is  a  polynomial in  $M$ whose degree will depend on the number of terms used in the power series approximation (which, in turn,  depends on the 
desired degree of accuracy of the approximation).

\subsection{Average hazard and reliability functions for PAR models}\label{AVhPAR}

Assuming a constant maintenance interval $M$, the time in which the equipment undertakes the $m$-maintenance activity is $mM$ and the age of the equipment in instant $t$ of period $m$, $w_{m-1}^+$, after the $(m-1)$-maintenance in equation (\ref{eq10}) becomes
\begin{equation}\label{eq10BIS}
\omega_{m-1}^+=t-\varepsilon (m-1)M. 
\end{equation}
The stationarity property that characterizes  the behaviour of age in the PAS model, does not hold for the PAR model. As a result  
the evaluation of the average hazard and reliability  for PAR models requires integration of the hazard and reliability functions over the replacement period $(0,RP)$.  Since the age representation given by \eqref{eq10BIS} is not continuous in $(0,RP)$ we use the following linear approximation, 
\begin{equation}\label{eq10TRIS}
\omega(t,\varepsilon)=t(1-\varepsilon)+\frac{M\varepsilon}{2},  \,\,\,\,\,0<t<RP.
\end{equation}
The age of the equipment (as a function of the chronological time) for the PAR model and the linear approximation in
\eqref{eq10TRIS} are shown in Figure 1-(b).

Substituting age by its expression given by \eqref{eq10TRIS}  in formulas \eqref{eq3} and  \eqref{eq7}, the hazard function for the PAR-linear and PAR-Weibull models becomes:
\begin{equation}\label{eqtt29}
h(t)=\left\{\begin{array}{ll}
\alpha\left(t(1-\varepsilon)+\frac{M\varepsilon}{2}\right),&\text{for PAR-linear},\\
\frac{\beta}{\eta^\beta}\left(t(1-\varepsilon)+\frac{M\varepsilon}{2}\right)^{\beta-1},&\text{for PAR-Weibull}.\end{array}\right.
\end{equation}
The  corresponding averaged hazard functions, $h^*$,  are obtained  integrating  \eqref{eqtt29} over  the replacement period  and dividing by $RP$,
\begin{equation}\label{eq34}
h^*(M)=
\frac{1}{RP}\int_{0}^{RP}h(t)\ dt=\left\{\begin{array}{ll}
\frac{\alpha}{2}(\varepsilon M+RP(1-\varepsilon)), &\text{ for PAR-linear},\\
\frac{(M\varepsilon+2RP(1-\varepsilon))^\beta-(M\varepsilon)^\beta}{RP(1-\varepsilon)(2\eta)^\beta}, &\text{ for PAR-Weibull}.\\
\end{array}\right.\\
\end{equation}
Similarly, substituting in (\ref{eqt39})  age by its expression in \eqref{eq10TRIS}, the  reliability function for the PAR-linear and PAR-Weibull models becomes
\begin{equation}\label{eqtt41}
R(t)=\left\{\begin{array}{ll}
\exp\left(-\frac{\alpha}{2}\left(t(1-\varepsilon)+\frac{M\varepsilon}{2}\right)^2\right), &\text{for PAR-linear},\\
\exp\left(-\left(\frac{t(1-\varepsilon)+\frac{M\varepsilon}{2}}{\eta}\right)^\beta\right),&\text{for PAR-Weibull}.
\end{array}\right.
\end{equation}
The average reliability function, $R^*$, for the PAR-linear and PAR-Weibull models is then obtained by integrating  
\eqref{eqtt41} over the interval $(0,RP)$ and dividing by the length of the replacement period, $RP$.
Once again, since the expressions in \eqref{eqtt41} are not integrable functions we can use their power series representation to perform the integration. 

\begin{figure}
    \centering
    \begin{subfigure}[t]{0.5\textwidth}
        \centering
        \includegraphics{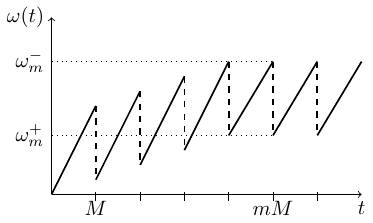}
        \caption{PAS.}
    \end{subfigure}%
    ~ 
    \begin{subfigure}[t]{0.5\textwidth}
        \centering
        \includegraphics{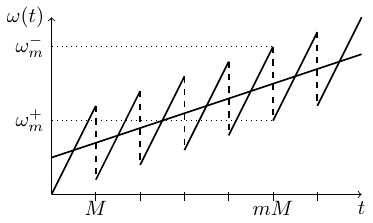}
        \caption{PAR.}
    \end{subfigure}
    \caption{PAS and PAR models age behaviour as function of time.}\label{fig1}
\end{figure}

\subsection{Cost models}

The relevant costs, in analyzing maintenance optimization of a safety-related equipment with a single component,  include the  cost contributions associated with performing preventive and corrective maintenance and the cost associated with replacing the component. In particular the hourly  cost of the equipment as a function of the maintenance interval $M$ is given by:
\begin{equation}\label{eq:C}
C(M)=\frac{c_m}{M} + \frac{1}{M}(\rho-h^*(M) M)c_c  + \frac{c_o}{RP},
\end{equation}
where $c_m$ and $c_c$ denote the cost associated with conducting each preventive and corrective maintenance, respectively;  $c_{o}$ is the total cost of replacing the component;   
$\rho$ is the cyclic or per-demand failure probability; $h^*$ is the average hazard function discussed in  subsections  \ref{AVhPAS}-\ref{AVhPAR};  and both $M$ and $RP$ are expressed in hours. 

The average reliability function, $R^*(\cdot)$,  and the hourly cost function, $C(\cdot)$,  that we derived here are used  as objective or restriction functions in  the optimization problem that we describe in more details in the next subsection.

\subsection{Multi-objective Optimization Problem}

The problem of optimizing the maintenance activities of an equipment has a natural formalization as \textit{Multi-Objective Optimization Problem} (MOP) where  the goal is 
to determine the value of the  maintenance intervals, over the replacement period, in each component of the equipment (the decision variables)  that maximize the average reliability and minimize the hourly cost (the objective functions) under a set of constraints, that, in our case, are restrictions on the maintenance intervals generated by  the value of the hourly cost and reliability, $(C_i , R_i)$, associated to the current maintenance intervals implemented in the plant. 
The solution of the MOP problem is a two step procedure. At the first step, the initial point $(C_i , R_i)$ determines the feasible region to solve  two \textit{Single-Objective Optimization Problems} (SOP), where only one criterion (cost or reliability) is involved in the optimization process acting as the single objective function, and the other  is implemented in the set of constraints. In particular,  the first SOP problem  consists in determine the value of the maintenance intervals in each component that minimize the cost per hour of the equipment 
subject to the restriction that the average reliability must be greater or equal to the initial value 
$R_{i}$. We denote the hourly cost and the average reliability associated to the solution of this first SOP problem  by $(C_r , R_o)$. The second SOP problem  consists in determine the value of the maintenance intervals in each component that maximize the average reliability  of the equipment 
subject to the restriction that the cost per hour must be smaller or equal to the initial value 
$C_{i}$. We denote the hourly cost and the average reliability associated to the solution of this second SOP problem  by $(C_o , R_r)$.

The solutions  of the two SOP's subproblems establish a new feasible region to solve the MOP acting as new restrictions and references to define the effectiveness of any feasible solution to the MOP  in terms of  the relative improvement in each objective function. We define the opposite of the relative improvement at each feasible solution as
\begin{equation*}\label{eq49}
e(C)=\frac{C-C_r}{C_r-C_o};\hspace{0.5cm}e(R)=\frac{R_r-R}{R_o-R_r}.
\end{equation*}

\begin{figure}
\begin{center}
\includegraphics[width=0.45\textwidth]{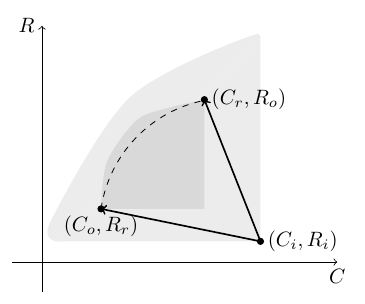}
\end{center}
\caption{Multi-Objective Optimization Problem.}\label{fig2}
\end{figure}

Applying the so-called weighted sum strategy, the multi-objective problem of minimizing the vector of objective functions is converted into a scalar problem by constructing a weighted sum of $e(C)$ and $e(R)$. In particular, if we take the weighting coefficients in interval $[0,1]$ wich sum to one, the  second step of the MOP procedure consists in the  minimization of a convex combination of both effectiveness functions, $e(C)$ and $e(R)$, subject to the restrictions generated by the solutions of each SOP.  The solution of this minimization problem produces the Pareto front of the MOP problem (see Figure \ref{fig2}) a rich set of non-dominated solutions from which the decision maker  can choose the one that provides  the best trade-off between cost and reliability.

The efficiency and accuracy of the solution of a  MOP problem depend, among other things, on the characteristics of the objective function and constraints.  In our case,  where  the objective function and constraints are nonlinear functions of the decision variables, the problem is known as \textit{Nonlinear Programming} (NLP). A solution of the NLP problem generally requires an iterative procedure to establish a direction search at each major iteration through the resolution of several \textit{Quadratic Programming}   subproblems (i.e. minimization of  a quadratic objective function that is linearly constrained). These methods are commonly referred to as \textit{Sequential Quadratic Programming} (SQP) methods. Details can be found in Biggs (1975), Powell (1983)
or Fletcher (1987).

\section{Application}\label{sec5}

In this section we illustrate the general approach to optimal maintenance activities described in this paper using a data set  containing failure time and 
maintenance activities of eight motor-operated safety valves of a Nuclear Power Plant (NPP) over a period of ten years. The failures are classified as electrical or mechanical. The electrical failures are associated with the actuator component while the mechanical failures correspond to the valve component. 
In order to mitigate the equipment aging process,  maintenance activities are currently  performed 
on each component every 180 days. The choice of 180 days for the maintenance interval of both components, however, is the result of heuristic considerations.

The goal is to use the observed data and the proposed approach to find  the optimal values of the maintenance intervals for each component  that possibly provide a better tradeoff between reliability and cost.

As described in the previous sections, for each component we consider four candidate models for the component behaviour: PAS-linear,  PAR-linear, PAS-Weibull,  and PAR-Weibull. Maximum likelihood  parameters estimates by model type and type of component  are shown in Table 1. 

In Table 1 we also show  the maximized likelihood ($L$) and  the values  of  the LCV,  the AIC, and BIC statistics used for model selection. Optimal values of the three selection criteria  are highlighted in bold.  As we noted in section  \ref{sec3:subsec2},  LCV is a positive oriented measure of model's performance (the larger the better) while AIC and BIC  are negatively oriented measures  (the smaller their value the better is the model). As shown in table 1 the results of model selection are robust with respect to the choice of the selection criterion. Regardless of the criterion used  (LCV, BIC or AIC) the best model for the actuator  is the PAS-Weibull, while the best model for the valve is the PAR-linear. Note that for the valve component  an analysis of model performance based on the maximized likelihood  $L$ (that ignores model complexity and over parametrization)  would select the PAR-Weibull model rather than the simpler PAR-linear model.
\begin{table}
\begin{center}
\begin{tabular}{c|c|c}
 & Actuator & Valve \\\hline
 PAS-linear & $\begin{array}{l}
\alpha=5.97\mathrm{e}{-9};\\ 
\varepsilon=1;\\
L = 2.03\mathrm{e}{-42}  ;\\  
LCV = 1.13\mathrm{e}{-39};\\
AIC = 195.9995   ;\\   
BIC = 199.3266
\end{array}$& $\begin{array}{l}
\alpha=1.54\mathrm{e}{-9};\\ \varepsilon=0.6002;\\
L = 1.57\mathrm{e}{-26}  ;\\  LCV = 7.98\mathrm{e}{-28};\\
AIC =122.8291    ;\\   BIC = 125.9398
\end{array}$\\\hline
PAR-linear & $\begin{array}{l}
\alpha=5.97\mathrm{e}{-9};\\  \varepsilon=1;\\
L = 2.03\mathrm{e}{-42}  ;\\   LCV = 1.13\mathrm{e}{-39};\\
AIC = 195.9995   ;\\   BIC =199.3266
\end{array}$& $\begin{array}{l}
\alpha=1.73\mathrm{e}{-9};\\ \varepsilon=0.7584;\\
L = 2.17\mathrm{e}{-26}  ;\\  {\bf LCV = 1.30\mathrm{e}{-27}};\\
{\bf AIC =122.1891}    ;\\   {\bf BIC =125.2998}\\
\end{array}$\\\hline
PAS-Weibull & $\begin{array}{l}
\beta=7.4708;\\  \eta=15397;\\ \varepsilon=0.8482;\\
{\bf L = 7.85\mathrm{e}{-40}}  ;\\  {\bf LCV = 9.09\mathrm{e}{-36}};\\
{\bf AIC = 186.0866}   ;\\   {\bf BIC =191.0773}
\end{array}$& $\begin{array}{l}
\beta=3.9505;\\  \eta=32465;\\ \varepsilon=0.3256;\\
L = 3.25\mathrm{e}{-26}  ;\\  LCV = 5.68\mathrm{e}{-30};\\
AIC = 123.3771   ;\\   BIC =128.0431
\end{array}$\\\hline
PAR-Weibull & $\begin{array}{l}
\beta=7.0403;\\  \eta=14971;\\  \varepsilon=0.9149;\\
L = 5.94\mathrm{e}{-40}  ;\\  LCV = 7.55\mathrm{e}{-36};\\
AIC =186.6421    ;\\   BIC =191.6328
\end{array}$& $\begin{array}{l}
\beta=3.8450;\\  \eta=29895;\\ \varepsilon=0.5282;\\
{\bf L = 4.47\mathrm{e}{-26}}  ;\\  LCV = 4.95\mathrm{e}{-30};\\
AIC = 122.7397   ;\\   BIC =127.4057
\end{array}$\\\hline
\end{tabular}
\end{center}
\caption{Summary of the estimation and model selection results}\label{tab1}
\end{table}
The best models for the actuator (PAS-Weibull) and for the valve (PAR-linear) are used as inputs for the optimization problem discussed in section \ref{sec4},
where, for both   components, optimization is performed under reliability and cost criteria, and maintenance intervals ($M_A$ for the actuator and $M_V$ for the valve)   are the  decision variables.  The  replacement period for both components is $RP = 87600$ hours.
The yearly cost  associated with conducting each preventive and corrective maintenance ($c_m$ and $c_c$),  the total cost of replacing the component ($c_{o}$) and the cyclic failure probability
 ($\rho$) by component type are shown in Table \ref{tab5}.
 \begin{table}
\begin{center}
\begin{tabular}{c|c|c|c|c}
 &  $\rho$ & $c_c$ &  $c_m$ & $c_o$ \\\hline
Actuator & $9.1\mathrm{e}{-4}$	& 3120	& 300	& 1900\\\hline
Valve & $9.1\mathrm{e}{-4}$	& 3120	& 800	& 3600\\\hline
\end{tabular}
\end{center}
\caption{Single cost data for actuator and valve in \euro/year.}\label{tab5}
\end{table}

Using the estimated  models in Table 1 (the PAS-Weibull for the actuator  and the PAR-linear  for the valve) and the cost parameters in Table 2, we derive  the average reliability and cost functions. 
As mentioned before, the average reliability function is obtained by integration of the corresponding reliability function by using a power series approximation. For the actuator  just two terms in the power series representation are sufficient to obtain an approximation $R_A^*(M_{A})$ of the average reliability function with an accuracy  of $10^{-9}$, 
\begin{equation*}\label{eq52}
R_A^*(M_{A})=1+\frac{\left(\frac{M_A}{\varepsilon\eta}\right)^\beta\left((1-\varepsilon)^{\beta+1}-1\right)}{\varepsilon_A(\beta+1)}.
\end{equation*}
For the   valve,  however, we  needed six terms in the power series representation of the reliability function to obtain an approximation $R_V^*(M_{V})$ with an accuracy of $10^{-7}$  (the expression of $R_V^*(M_{V})$, a polynomial of  degree 10 in $M_V$ whose coefficients are products of powers of parameters,  is omitted due to its extension and complexity).
The global average reliability function of the equipment, $R^*(M_{A},M_{V})$, is then obtained  as the product of the average reliability function  of each component, 
\begin{equation}\label{R:star}
R^*(M_{A},M_{V}) = R_A^* (M_{A}) \cdot R_V^*(M_{V}).
\end{equation}

Replacing  in \eqref{eq:C}, the cost parameters in Table 2 and the expression of the averaged hazard functions, $h^*$, in  \eqref{eq34} for the selected models  (PAS-Weibull for the actuator  and the PAR-linear  for the valve), we obtain the expression of hourly cost of the equipment as a function of the decision variables $M_{A}$ and   $M_{V}$ (to unify time units to hours we multiply yearly cost by  $8760$), 
\begin{equation}\label{eq53}
\begin{array}{ll}
C(M_{A},M_{V})=& 8760\left[\frac{1}{M_A}\left(300+3120\left(9.1\mathrm{e}{-4}+\left(\frac{M_A}{\varepsilon_A\eta}\right)^\beta(1-(1-\varepsilon)^\beta)\right)\right)\right. \nonumber \\ 
& + \left.\frac{1}{M_V}\left(800+3120\left(9.1\mathrm{e}{-4}+\frac{\alpha M_V}{2}(\varepsilon_V M_V+L(1-\varepsilon_V)\right)\right)\right]\\
&+550.
\end{array}
\end{equation}

Finally, using (\ref{R:star}), (\ref{eq53})  and the current  maintenance activities intervals implemented in the NPP  for the actuator and the valve ($M_A = M_V =$4320 hours, i.e. 180 days), 
we obtain  the initial values for the objective functions:  $R_i = 0.857848$ and $C_i = 3372.94$\euro. These values act as restrictions to solve the single-objective optimization subproblems, which are stated as:
\begin{equation}\label{eq54}
\begin{array}{ll}
\min & C(M_A,M_V)\\
\text{s.t.} & R\geq R_i\\
& 0\leq M_A,M_V\leq 87600
\end{array} 
\text{ and }
\begin{array}{ll}
\min & -R(M_A,M_V)\\
\text{s.t.} & C\leq C_i\\
& 0\leq M_A,M_V\leq 87600.
\end{array} 
\end{equation}

The solution of the cost minimization subproblem (in days) is  $(M_{A},M_{V})=(270,176)$, and  
reduces the yearly cost in nearly 150\euro\ in each equipment, keeping their averaged reliability slightly higher than its initial value. 
The solution of the reliability maximization subproblem (in days) is  $(M_{A},M_{V})=(261,162)$ with an  associated reliability of 0.860161 and an associated cost lower than the current one. The solutions of both subproblems generate the new feasible set where MOP must be solved. Taking values of weighting coefficients in the interval [0,1] and decision variables bounded in the interval [0,87600], the SQP method provides 155 non dominated solutions which describe the Pareto front of our problem. Table \ref{tab6} shows some of these solutions in addition to those associated to the  initial/current values, $(M_{A},M_{V})=(180,180)$,  and to the solutions of the subproblems in (\ref{eq54}).

\begin{table}
\begin{center}
\begin{tabular}{c|c|c|c|c}
 &  $M_A$ & $M_C$ &  $C$ & $R$ \\\hline
Initial & 180	& 180	& 3372.94	& 0.8578\\\hline
& $\dots$	& $\dots$	& $\dots$	& $\dots$\\\hline
$opt(R)$ & 261	&162	&3371.89	&0.8602\\\hline
& $\dots$	& $\dots$	& $\dots$	& $\dots$\\\cline{2-5}
&264&	165	&3336.87&	0.8597\\\cline{2-5}
& $\dots$	& $\dots$	& $\dots$	& $\dots$\\\cline{2-5}
&265&	167&	3229.51&	0.8579\\\cline{2-5}
& $\dots$	& $\dots$	& $\dots$	& $\dots$\\\cline{2-5}
&266&	169&	3294.38&	0.8590\\\cline{2-5}
& $\dots$	& $\dots$	& $\dots$	& $\dots$\\\cline{2-5}
&267&	172&	3264.38&	0.8585\\\cline{2-5}
& $\dots$	& $\dots$	& $\dots$	& $\dots$\\\cline{2-5}
&268&	174&	3244.62&	0.8582\\\cline{2-5}
& $\dots$	& $\dots$	& $\dots$	& $\dots$\\\hline
$opt(C)$ & 270	&176	&3224.35	&0.8579\\\hline
\end{tabular}
\end{center}
\caption{Some results of the optimization process.}\label{tab6}
\end{table}

\section{Conclusions}\label{sec6}

In this paper we describe  a general approach to optimal imperfect maintenance activities of a repairable equipment with $q$ independent components.  The proposed procedure allows to relax  some of the assumptions that characterize existing proposals. In particular rather than assuming that  all the components of the equipment share the same failure rate and maintenance models and the same maintenance interval, and that the effectiveness parameter is known, we consider for each component  a class of four candidate models that are obtained combining two possible models for failure rate (Weibull or linear) and two possible models for imperfect maintenance (PAS y PAR) and let the data select the most suitable model for each component.  Both the parameters of the failure rate model and  the effectiveness parameter of the IM models are assumed to be unknown and are jointly estimated via maximum likelihood. Model selection is performed, separately for each component, using three standard criteria that take into account the problem of over-parametrization.  
The selected models are used  to derive the cost per unit time and the average reliability of the equipment, the objective functions of a Multi-Objective Optimization Problem with maintenance intervals of the single components as decision variables.

We illustrate the advantages of the proposed procedure using a data set  containing failure time and maintenance activities of eight motor-operated safety valves of a NPP over a period of ten years. The results of this real data example well illustrate: (1)  the importance of considering a class of model for each component and allow the selected model and the maintenance interval to be different for different components of an equipment (the best model for actuator and valve  in the example are 
PAS-Weibull and  PAR-linear respectively, and the optimal maintenance intervals for the two components that make the Pareto front of the optimization problem are also different);  (2) the limitations of the maximized likelihood as model selection criterion
and the importance of considering alternative criteria that  do take into account model complexity and over-parametrization (in the example the maximized likelihood would pick the Weibull-PAR model as best model for the valve rather than the simpler linear-PAR model selected by AIC, BIC and LVC criteria); (3) the importance of assuming the effectiveness parameter of the IM models to be unknown (notice how the estimates of $\varepsilon$ in Table~\ref{tab1} differ across  the different models of the two components).

\section{Acknowledgments}
This document is a collaborative effort and is partially supported by Conselleria d'Educaci\'o, Investigaci\'o, Cultura i Sport (Generalitat de la Comunitat Valenciana, Spain) under grant GV/2017/015.


\section{References}


\begin{thebibliography}{99}

\bibitem{A74} Akaike, H., 1974. A new look at the statistical model identication. IEEE Transactions on Automatic Control 19, 716--723.
\bibitem{AT15} Alrabghi, A., Tiwari, A., 2015. State of the art in simulation-based optimization for maintenance systems. Computers \& Industrial Engineering
82, 167--182.
\bibitem{B75} Biggs, M., 1975. Constrained minimization using recursive quadratic programming. Towards Global Optimization, 341--349.
\bibitem{FBCR17} Feng, Q., Bi, W., Chen, Y., Ren, Y., 2017. Cooperative game approach based on agent learning for  eet maintenance oriented to mission reliability. Computers \& Industrial Engineering 112, 221--230.
\bibitem{F87} Fletcher, R., 1987. Practical Methods for Optimization. John Willey and Sons, London.
\bibitem{HTF08} Hastie, T., Tibshirani, R., Friedman, J., 2008. The Elements of Statistical Learning. Springer, New York.

\bibitem{K89} Kijima, M., 1989. Some results for repairable systems with general repair. Journal of Applied Probability 26, 89--102.
\bibitem{KS86} Kijima, M., Sumita, N., 1986. A useful generalization of renewal theory: counting process governed by non-negative markovian increments. Journal
of Applied Probability 23, 71--88.
\bibitem{LRWW98} Lagarias, J., Reeds, J., Wright, M., Wright, P., 1998. Convergence properties of the nelder-mead simplex method in low dimensions. SIAM Journal of
Optimization 9, 112--147.
\bibitem{M79} Malik, M., 1979. Reliable preventive maintenance scheduling. AIIE Transactions 11, 221--228.

\bibitem{MSS99} Martorell, S., S\'anchez, A., Serradell, V., 1999. Age-dependent reliability model considering effects of maintenance and working conditions. Reliability Engineering and System Safety 64, 19--31.
\bibitem{NM65} Nelder, J., Mead, R., 1965. Simplex method for function minimization. Computer Journal, 308--313.
\bibitem{NSM11} Neves, M., Santiago, L., Maia, C., 2011. A condition-based maintenance policy and input parameters estimation for deteriorating systems under
periodic inspection. Computer \& Industrial Engineering 61, 503--511.
\bibitem{P83} Powell, M., 1983. Variable metric methods for constrained optimization. Mathematical Programming. The State of the Art, 288--311.


\bibitem{QCG17} Qiu, Q., Cui, L., Gao, H., 2017a. Availability and maintenance modelling for systems subject to multiple failure modes. Computers \& Industrial
Engineering 108, 192--198.

\bibitem{QCSY17} Qiu, Q., Cui, L., Shen, J., Yan, L., 2017b. Optimal maintenance policy considering maintenance errors for systems operating under performance-
based contracts. Computers \& Industrial Engineering 112, 147--155.

\bibitem{S78} Schwarz, G., 1978. Estimating the dimension of a model. Annals of Statistics 6, 461--464.

\bibitem{TB13} Taghipour, S., Banjevic, D., 2013. Maximum likelihood estimation from interval censored recurrent event data. Computers \& Industrial Engineering
64, 143--152.

\bibitem{TK95} Tan, J., Kramer, M., 1995. Hazard function modelling using cross validation: from data collection to model selection. Reliability Engineering and System Safety 49, 155--169.

\bibitem{TRB14} Tanwar, M., Rai, R., Bolia, N., 2014. Imperfect repair modeling using kijima type generalized renewal process. Reliability Engineering and System
Safety 124, 24--31.

\bibitem{TLL11} Tsai, T.-R., Liu, P.-H., Lio, Y., 2011. Optimal maintenance time for imperfect maintenance actions on repairable product. Computers \& Industrial
Engineering 60, 744--749.

\bibitem{WP06} Wang, H., Pham, H., 2006. Reliability and Optimal Maintenance. Springer-Verlag, London.

\bibitem{WY12} Wang, Z.-M., Yang, J.-G., 2012. Numerical method for weibull generalized renewal process and its applications in reliability analysis of nc machine tools. Computers \& Industrial Engineering 63, 1128--1134.

\bibitem{XXZD12} Xia, T., Xi, L., Zhou, X., Du, S., 2012. Modeling and optimizing maintenance schedule for energy subject to degradation. Computers \& Industrial
Engineering 63, 607--614.

\bibitem{ZZZW16} Zheng, Z., Zhou, W., Zheng, Y., Wu, Y., 2016. Optimal maintenance policy for a system with preventive repair and two types of failures. Computers
\& Industrial Engineering 98, 102--112.
\end{thebibliography}
\end{document}